\documentclass{amsart}

\usepackage{fancyhdr}
\usepackage{lastpage}
\usepackage{stmaryrd,yhmath}

\newcommand{\bdis}{\begin{displaymath}}
\newcommand{\edis}{\end{displaymath}}
\newcommand{\be}{\begin{equation}}
\newcommand{\ee}{\end{equation}}
\newcommand{\mbb}{\mathbb}
\newcommand{\mcal}{\mathcal}

\newcommand{\vp}{\varphi}

\newcommand{\zf}{\zeta\left(\frac{1}{2}+it\right)}

\theoremstyle{definition}

\theoremstyle{remark}
\newtheorem{remark}[]{Remark}

\newtheorem*{mydef1}{{\bf Theorem}}

\newtheorem*{mydef2}{{\bf Definition}}

\newtheorem*{mydef5}{{\bf Lemma}}

\numberwithin{equation}{section}

\begin{document}

\title{Jacob's ladders and the three-points interaction of the Riemann zeta-function with itself}

\author{Jan Moser}

\address{Department of Mathematical Analysis and Numerical Mathematics, Comenius University, Mlynska Dolina M105, 842 48 Bratislava, SLOVAKIA}

\email{jan.mozer@fmph.uniba.sk}


\begin{abstract}
It is proved that some set of the values of $|\zeta(\sigma_0+i\vp_1(t))|^2$  on every fixed line $\sigma=\sigma_0>1$ generates a corresponding
set of the values of $|\zeta(\frac 12+it)|^2$ on the critical line $\sigma=\frac 12$ (i.e. we have an analogue of the Faraday law).
\end{abstract}

\maketitle

\section{The result}

\subsection{}

Let us remind the theorem of H. Bohr (see \cite{5}, p. 253, Theorem 11.6(C)): The function $\zeta(s)$ attains every value $a\in\mbb{C}$ except
$0$ infinitely many times in the strip $a<\sigma<1+\delta$. At the same time, this theorem does not allow us to prove anything about the set of the
roots of the equation
\bdis
|\zeta(\sigma_0+it)|=|a|,\ \sigma_0\in (1,1+\delta) ,
\edis
i.e. about the roots lying on every fixed line $\sigma=\sigma_0$. \\

In this paper we will study more complicated nonlinear equation

\be \label{1.1}
\left|\zeta\left(\frac 12+iu\right)\right|^2\left|\zeta\left(\frac 12+iv\right)\right|^2=\frac 12\zeta(2\sigma)\ln u,\ u,v>0,\ \sigma\geq\alpha>1
\ee
where $\alpha$ is an arbitrary fixed value.

\begin{remark}
The theory of H. Bohr is not applicable on the equation (\ref{1.1}). However, by the method of Jacob's ladders we obtain some information about the set
of approximative solutions of this equation.
\end{remark}

\subsection{}

Let us remind that

\be \label{1.2}
\tilde{Z}^2(t)=\frac{{\rm d}\vp_1(t)}{{\rm d}t},\ \vp_1(t)=\frac 12\vp(t)
\ee

where

\be \label{1.3}
\tilde{Z}^2(t)=\frac{Z^2(t)}{2\Phi'_\vp[\vp(t)]}=\frac{Z^2(t)}{\left\{ 1+\mcal{O}\left(\frac{\ln\ln t}{\ln t}\right)\right\}\ln t} ,
\ee
(see \cite{1}, (3.9); \cite{2}, (1.3); \cite{4}, (1.1), (3.1)), and $\vp(t)$ is the solution of the nonlinear integral equation

\bdis
\int_0^{\mu[x(T)]}Z^2(t)e^{-\frac{2}{x(T)}t}{\rm d}t=\int_0^TZ^2(t){\rm d}t .
\edis

\subsection{}

\begin{mydef2}
If there are some sequences $\{ u_n(\sigma_0)\}_{n=0}^\infty,\ \{ v_n(\sigma_0)\}_{n=0}^\infty,\ \sigma_0>1$ for which the following condition
are fulfilled
\be \label{eqA}
\lim_{n\to\infty} u_n(\sigma_0)=\infty,\ \lim_{n\to\infty}v_n=\infty , \tag{A}
\ee
\be \label{eqB}
\frac{\left|\zeta\left(\frac 12+iu_n(\sigma_0)\right)\right|^2\left|\zeta\left(\sigma_0+iv_n(\sigma_0)\right)\right|^2}
{\zeta(2\sigma_0)\ln u_n(\sigma_0)}=\frac 12+o(1) , \tag{B}
\ee
then we call each ordered pair $[u_n(\sigma_0),v_n(\sigma_0)]$ an \emph{asymptotically approximate solution} (AA solution) of the nonlinear
equation (\ref{1.1}).
\end{mydef2}

The following theorem holds true.

\begin{mydef1}
Let us define the continuum set of sequences
\bdis
\{ K_n(T)\}_{n=0}^\infty,\ K_n\geq T_0[\vp_1]
\edis
as follows
\be \label{1.4}
\begin{split}
& K_0=T,\ K_1=K_0+K_0^{1/3+2\epsilon},\ K_2=K_1+K_1^{1/3+2\epsilon}, \dots , \\
& K_{n+1}=K_n+K_n^{1/3+2\epsilon} , \dots \ .
\end{split}
\ee
Then for every $\sigma_0:\ \sigma_0\geq \alpha>1$ there is a sequence $\{ u_n(\sigma_0)\}_{n=0}^\infty,\ u_n(\sigma_0)\in (K_n,K_{n+1})$ such that
\be \label{1.5}
\frac{\left|\zeta\left(\frac 12+iu_n(\sigma_0)\right)\right|^2\left|\zeta\left(\sigma_0+i\vp_1[u_n(\sigma_0)]\right)\right|^2}
{\zeta(2\sigma_0)\ln u_n(\sigma_0)}=\frac 12+\mcal{O}\left(\frac{\ln\ln K_n}{\ln K_n}\right) .
\ee
holds true where $\vp_1[u_n(\sigma_0)]\in (\vp_1(K_n),\vp_1(K_{n+1}))$ and
\be \label{1.6}
\rho\{[K_n,K_{n+1}];[\vp_1(K_n),\vp_1(K_{n+1})]\}\sim (1-c)\pi(K_n)\to\infty
\ee
as $n\to\infty$ where $\rho$ denotes the distance of the segments, $c$ is the Euler's constant and $\pi(t)$ is the prime-counting function, i.e. the
ordered pair
\bdis
[u_n(\sigma_0),\vp_1[u_n(\sigma_0)]]; \ \vp_1[u_n(\sigma_0)]=v_n(\sigma_0)
\edis
is the AA solution of the nonlinear equation (\ref{1.1}).
\end{mydef1}

\begin{remark}
Let us point out that the formula (\ref{1.5}) binds together the values of $\zeta(s)$ at three distinct points: at the point $\frac 12+iu_n(\sigma_0)$
lying on the critical line $\sigma=\frac 12$, and at the points $\sigma_0+i\vp_1[u_n(\sigma_0)],2\sigma_0$ lying in the semi-plane
$\sigma\geq\alpha>1$ .
\end{remark}

\begin{remark}
Since by the eq. (\ref{1.5}) we have
\be \label{1.6}
\left|\zeta\left(\frac 12+iu_n(\sigma_0)\right)\right|^2\sim \frac{\zeta(2\sigma_0)\ln u_n(\sigma_0)}{2}\frac{1}
{\left|\zeta\left(\sigma_0+i\vp_1[u_n(\sigma_0)]\right)\right|^2} ,
\ee
\be \label{1.7}
\left|\zeta\left(\sigma_0+i\vp_1[u_n(\sigma_0)]\right)\right|^2\sim \frac{\zeta(2\sigma_0)\ln u_n(\sigma_0)}{2}\frac{1}
{\left|\zeta\left(\frac 12+iu_n(\sigma_0)\right)\right|^2} ,
\ee
then for the two parallel conductors placed at the positions of the lines $\sigma=\frac 12$ and $\sigma=\sigma_0\geq \alpha>1$ some art of the
Faraday law holds true: the sequence of the energies $\{|\zeta(\sigma_0+i\vp_1[u_n(\sigma_0)])|^2\}_{n=0}^\infty$ on $\sigma=\sigma_0$ generates
the sequence of the energies $\{|\zeta(\frac 12+iu_n(\sigma_0))|^2\}_{n=0}^\infty$ on $\sigma=\frac 12$, and vice versa (see (\ref{1.6}), (\ref{1.7});
the energy is proportional to the square of the amplitude of oscillations).
\end{remark}

\section{The local mean-value theorem}

\subsection{}

Let us remind that there is the \emph{global} mean-value theorem (see \cite{5}, p. 116)

\bdis
\lim_{T\to\infty}\frac{1}{2T}\int_{-T}^T |\zeta(\sigma+it)|^2{\rm d}t=\zeta(2\sigma),\ \sigma>1 .
\edis

However, for our purpose, we need the \emph{local} mean-value theorem, i.e. the formula for the the integral

\bdis
\int_T^{T+U}|\zeta(\sigma+it)|^2{\rm d}t ,\ \sigma>1 .
\edis

In this direction, the following lemma holds true.

\begin{mydef5}
The formula
\be \label{2.1}
\int_T^{T+U}|\zeta(\sigma+it)|^2{\rm d}t=\zeta(2\sigma)U+\mcal{O}(1)
\ee
holds true uniformly for $T,U>0,\ \sigma\geq\alpha$, where $\alpha>1$ is an arbitrary fixed value. The $\mcal{O}$-constant depends of course on the choice
of $\alpha$.
\end{mydef5}

\begin{remark}
The formula (\ref{2.1}) is the asymptotic formula for $U\geq \ln\ln T$, for example.
\end{remark}

\subsection{Proof of the Lemma}

Following the formula

\bdis
\zeta(s)=\zeta(\sigma+it)=\sum_{n=1}^\infty \frac{1}{n^{\sigma+it}}, \ \sigma>1
\edis

we obtain

\bdis
|\zeta(\sigma+it)|^2=\zeta(2\sigma)+\sum_n\sum_{m\not=n}\frac{1}{(mn)^\sigma}\cos\left( t\ln\frac nm\right) ,
\edis
and
\be \label{2.2}
\begin{split}
& \int_T^{T+U}|\zeta(\sigma+it)|^2{\rm d}t=\zeta(2\sigma)U+\mcal{O}\left(\sum_{n}\sum_{m<n}\frac{1}{(mn)^\sigma\ln\frac nm}\right)= \\
& = \zeta(2\sigma)U+S(m<n)
\end{split}
\ee
uniformly for $T,U>0$. Let

\be \label{2.3}
S(m<n)=S\left( m<\frac n2\right)+S\left( \frac n2\leq m<n\right)=S_1+S_2 .
\ee

Since we have $2<\frac nm$ in $S_1$, and
\bdis
\sum_{n=1}^\infty\frac{1}{n^\sigma}=1+\sum_{n=2}^\infty \frac{1}{n^\sigma}<1+\int_1^\infty x^{-\sigma}{\rm d}x=1+\frac{1}{\sigma-1}
\edis
then
\be \label{2.4}
S_1=\mcal{O}\left(\sum_n\sum_{m<\frac n2}\frac{1}{(mn)^\sigma}\right)=\mcal{O}\left\{\left(\sum_{n=1}^\infty\frac{1}{n^\sigma}\right)^2\right\}=
\mcal{O}(1) .
\ee

We put $m=n-r>\frac n2;\ 1\leq r<\frac n2$ in $S_2$ and we obtain as usual
\bdis
\ln\frac nm=\ln\frac{n}{n-r}=-\ln\left( 1-\frac rn\right)>\frac rn .
\edis
Next, we have
\bdis
\begin{split}
& \sum_{n=2}^\infty \frac{\ln n}{n^{2\sigma-1}}=\frac{\ln 2}{2^{2\sigma-1}}+\sum_{n=3}^\infty \frac{\ln n}{n^{2\sigma-1}}<
\frac{\ln 2}{2^{2\sigma-1}}+\int_2^\infty\frac{\ln x}{x^{2\sigma-1}}{\rm d}x= \\
& =\frac{\ln 2}{2^{2\sigma-1}}+\frac{2^{-2\sigma+1}}{\sigma-1}\ln 2+\frac{2^{-2\sigma}}{(\sigma-1)^2}=\mcal{O}(1);\ 2^{-\sigma}\in \left( 0,\frac 12\right) .
\end{split}
\edis
Then
\be \label{2.5}
S_2=\mcal{O}\left(\sum_n\sum_{r=1}^{n/2}\frac{n}{(mn)^\sigma r}\right)=\mcal{O}\left( 2^\sigma\sum_{n=2}^\infty \frac{\ln n}{n^{2\sigma-1}}\right)=\mcal{O}(1) .
\ee
Finally, from (\ref{2.2}) by (\ref{2.3})-(\ref{2.5}) the formula (\ref{2.1}) follows.

\section{Proof of the Theorem}

\subsection{}

Let us remind that the following lemma holds true (see \cite{3}, (2.5);\cite{4}, (3.3)): for every integrable function (in the Lebesgue sense)
$f(x),\ x\in [\vp_1(T),\vp_1(T+U)]$ the following is true
\be \label{3.1}
\int_T^{T+U}f[\vp_1(t)]\tilde{Z}^2(t){\rm d}t=\int_{\vp_1(T)}^{\vp_1(T+U)}f(x){\rm d}x ,\ U\in \left(\left. 0,\frac{T}{\ln T}\right]\right. ,
\ee
where $t-\vp_1(t)\sim (1-c)\pi(t)$.

\subsection{}

In the case $f(t)=|\zeta(\sigma_0+it)|^2,\ U=U_0=T^{1/3+2\epsilon}$ we obtain from (\ref{3.1}) the following formula

\be \label{3.2}
\int_T^{T+U_0}|\zeta(\sigma_0+i\vp_1(t))|^2\tilde{Z}^2(t){\rm d}t=\int_{\vp_1(T)}^{\vp_1(T+U_0)}|\zeta(\sigma_0+it)|^2{\rm d}t .
\ee
Since (see (\ref{2.1}))
\be \label{3.3}
\begin{split}
& \int_{\vp_1(T)}^{\vp_1(T+U_0)}|\zeta(\sigma_0+it)|^2{\rm d}t=\zeta(2\sigma_0)\{\vp_1(T+U_0)-\vp_1(T)\}+\mcal{O}(1)= \\
& =\frac 12\zeta(2\sigma_0)U_0\tan[\alpha(T,U_0)]+\mcal{O}(1)
\end{split}
\ee
where (see (\ref{1.2}))
\bdis
\frac{\vp_1(T+U_0)-\vp_1(T)}{U_0}=\frac 12\frac{\vp(T+U_0)-\vp(T)}{U_0}=\frac 12\tan[\alpha(T,U_0)] ,
\edis
and (see \cite{3}, (2.6))
\bdis
\tan[\alpha(T,U_0)]=1+\mcal{O}\left(\frac{1}{\ln T}\right) ,
\edis
then (see (\ref{3.3}))
\be \label{3.4}
\int_{\vp_1(T)}^{\vp_1(T+U_0)}|\zeta(\sigma_0+it)|^2{\rm d}t=\frac 12\zeta(2\sigma_0)U_0\left\{ 1+\mcal{O}\left(\frac{1}{\ln T}\right)\right\} .
\ee

\subsection{}

Next, by the first application of the mean-value theorem we obtain (see (\ref{1.3}), (\ref{3.2}))

\be \label{3.5}
\begin{split}
& \int_T^{T+U_0}|\zeta(\sigma_0+i\vp_1(t))|^2\tilde{Z}^2(t){\rm d}t= \\
& =\frac{1}{\left\{ 1+\mcal{O}\left(\frac{\ln\ln\xi_1}{\ln \xi_1}\right)\right\}\ln\xi_1}
\int_T^{T+U_0}|\zeta(\sigma_0+i\vp_1(t))|^2\left|\zf\right|^2{\rm d}t , \\
& \xi_1=\xi_1(\sigma_0;T,U_0)\in (T,T+U_0) ,
\end{split}
\ee
and by the second application of this we obtain

\be \label{3.6}
\begin{split}
& \int_T^{T+U_0}|\zeta(\sigma_0+i\vp_1(t))|^2\left|\zf\right|^2{\rm d}t= \\
& =|\zeta(\sigma_0+i\vp_1(\xi_2))|^2
\left|\zeta\left(\frac 12+i\xi_2\right)\right|^2U_0 , \\
& \xi_2=\xi_2(\sigma_0;T,U_0)\in (T,T+U_0),\ \vp_1(\xi_2)\in (\vp_1(T),\vp_1(T+U));\\
& \ln\xi_1\sim\ln\xi_2 .
\end{split}
\ee

Hence, from (\ref{3.5}) by (\ref{3.6}) we have

\be \label{3.7}
\int_T^{T+U_0}|\zeta(\sigma_0+i\vp_1(t))|^2\tilde{Z}^2(t){\rm d}t=
\frac{|\zeta(\sigma_0+i\vp_1(\xi_2))|^2\left|\zeta\left(\frac 12+i\xi_2\right)\right|^2}
{\left\{ 1+\mcal{O}\left(\frac{\ln\ln\xi_2}{\ln \xi_2}\right)\right\}\ln\xi_2}U_0 ,
\ee
and from (\ref{3.2}) by (\ref{3.4}), (\ref{3.7}) we obtain

\be \label{3.8}
\frac{|\zeta(\sigma_0+i\vp_1(\xi_2))|^2\left|\zeta\left(\frac 12+i\xi_2\right)\right|^2}{\zeta(2\sigma_0)\ln\xi_2}=
\frac 12+\mcal{O}\left(\frac{\ln\ln\xi_2}{\ln\xi_2}\right) .
\ee

\subsection{}

Now, if we apply (\ref{3.8}) in the case (see (\ref{1.4}))

\bdis
\begin{split}
& [T,T+U_0]\to [K_n,K_{n+1}];\ \xi_2(\sigma_0)\to\xi_{2,n}(\sigma_0)\in (K_n,K_{n+1}); \\
& \xi_{2,n}(\sigma_0)=u_n(\sigma_0) ,
\end{split}
\edis
then we obtain (\ref{1.5}). The statement (\ref{1.6}) follows from $t-\vp_1(t)\sim (1-c)\pi(t)$.

\section{Concluding remarks}

Let us remind that in \cite{1} we have shown the following formula holds true

\bdis
\int_0^T Z^2(t){\rm d}t=\vp_1(T)\ln\vp_1(T)+(c-\ln 2\pi)\vp_1(T)+c_0+\mcal{O}\left(\frac{\ln T}{T}\right),\ \vp_1(T)=\frac{\vp(T)}{2} ,
\edis
where $\vp_1(T)$ is the Jacob's ladder. It is clear that $\vp_1(T)$ is the asymptotic solution of the nonlinear transcendental equation
\bdis
\int_0^T Z^2(t){\rm d}t=V(T)\ln V(T)+(c-\ln 2\pi)V(T) .
\edis


\begin{thebibliography}{29}% Replace 9 by 99 if 10 or more references
%
% Please note the use of "\and" between author names below
%\bibitem{fock}
%V.A. Fock, `The theory of space, time and gravitation`, GITTL, Moscow, 1955, (in russian).
%
\bibitem{1}
J. Moser, `Jacob's ladders and the almost exact asymptotic representation of the Hardy-Littlewood integral', Math. Notes 2010, {\bf 88},
pp. 414-422, arXiv: 0901.3973.
%
\bibitem{2}
J. Moser, `Jacob's ladders and the multiplicative asymptotic formula for short and microscopic parts of the Hardy-Littlewood integral', (2009),
arXiv:0907.0301.
%
\bibitem{3}
J. Moser,
`Jacob's ladders and the first asymptotic formula for the expression of the fifth order
$Z[\varphi(t)/2+\rho_1]Z[\varphi(t)/2+\rho_2]Z[\varphi(t)/2+\rho_3]\hat{Z}^2(t)$ for the collection of disconnected sets`, (2009),
arXiv:0912.0130.
%
\bibitem{4}
J. Moser, `Jacob's ladders, the iterations of Jacob's ladder $\vp_1^k(t)$ and asymptotic formulae for the integrals of the products
$Z^2[\varphi^n_1(t)]Z^2[\varphi^{n-1}(t)]\cdots Z^2[\varphi^0_1(t)]$ for arbitrary fixed $n\in \mbb{N}$` (2010), arXiv:1001.1632.
%
\bibitem{5}
E.C. Titchmarsh, `The theory of the Riemann zeta-function` Clarendon Press, Oxford, 1951.
%

\end{thebibliography}
\end{document}